\documentclass[12pt]{amsart}
\usepackage{amssymb,latexsym,comment,url}
\usepackage{amsfonts}
\usepackage{amssymb}
\usepackage{amsthm}
\usepackage{newlfont}
\usepackage{graphicx}
\usepackage{t1enc}
\usepackage{cite}
\usepackage[latin1]{inputenc}
\usepackage[english]{babel}

\newcommand{\bq }{\begin{equation}}
\newcommand{\eq }{\end{equation}}

\newenvironment{Proof}{\par\noindent{\sc Proof:}}%
                      {\hspace*{\fill}\nobreak$\Box$\par\medskip}
                       {\hspace*{\fill}\nobreak$\Box$\par\medskip}

\newtheorem{Proposition}{Proposition}[section]
\newtheorem{Theorem}[Proposition]{Theorem}
\newtheorem{Lemma}[Proposition]{Lemma}
\newtheorem{Corollary}[Proposition]{Corollary}

\theoremstyle{definition}

\begin{document}

\title[Fibonacci and Chebyshev Polynomials]%
{New Hypergeometric Connection Formulae Between Fibonacci and
Chebyshev Polynomials}

\author[W. M. Abd-Elhameed]%
{Waleed~M. Abd-Elhameed}
\address{Department of Mathematics, Faculty of Science,
Cairo University, Giza, Egypt}
\email{walee$_-$9@yahoo.com}

\author[Y. H. Youssri]%
{Youssri~H. Youssri}
\address{Department of Mathematics, Faculty of Science,
Cairo University, Giza, Egypt}
\email{youssri@sci.cu.edu.eg}

\author[N. El-Sissi]%
{Nermine~El-Sissi}
\address{American University in Cairo, Mathematics and Actuarial Science Department, AUC Avenue, New Cairo, Egypt}
\email{nelsissi@aucegypt.edu}

\author[M. Sadek]%
{Mohammad~Sadek}
\address{American University in Cairo, Mathematics and Actuarial Science Department, AUC Avenue, New Cairo, Egypt}
\email{mmsadek@aucegypt.edu}

\date{4th April 2015}
\maketitle

\begin{abstract}
  We establish new connection formulae between Fibonacci polynomials and
Chebyshev polynomials of the first and second kinds. These formulae
are expressed in terms of certain values of hypergeometric functions of the type
$_2F_{1}$. Consequently, we obtain some new expressions for the
celebrated Fibonacci numbers and their
 derivatives sequences. Moreover, we evaluate some definite integrals involving products of Fibonacci and
Chebyshev polynomials.
\end{abstract}

\let\thefootnote\relax

\section{Introduction}

Fibonacci numbers appear in several disciplines of modern science. The wide spectrum of applications of these numbers in mathematics, computer science, physics, biology, graph theory and statistics justifies the growing interest of mathematicians in the properties enjoyed by these numbers. The beautiful book of Koshy, \cite{koshy2011fibonacci}, exhibits some of the applications in which these numbers arise. 

The family of Fibonacci polynomials $\displaystyle\{F_n(x)\}$ is defined by Fibonacci-like recurrence relations. In fact, the sequence of Fibonacci numbers can be obtained from the sequence of Fibonacci polynomials by setting $x=1$. Therefore, the more knowledge we acquire on Fibonacci polynomials, the closer we get to understanding the qualities of Fibonacci numbers and other sequences of numbers. Yet, studying Fibonacci polynomials for their own sake provides us with a clearer idea concerning their combinatorial and analytic properties. A great deal of mathematical ingenuity has been invested in developing identities involving Fibonacci polynomials, Fibonacci numbers, and their generalizations, see
\cite{gulec2013new,taskara2010properties,
 el2010new, falcon2007fibonacci,yazlik2012note,wloch2013some} and the references there for examples on such identities.
 
 In order to study Fibonacci polynomials, one may consider linking Fibonacci polynomials to other well-studied polynomials, such as Chebyshev polynomials. Chebyshev
polynomials of the first and second kinds, $T_{n}(x)$ and $U_{n}(x)$ respectively, are subfamilies of the larger class of Jacobi polynomials. They are of crucial importance from both the theoretical and practical points of view. The interested
reader in these polynomials may consult \cite{mason2010chebyshev}.

 Given two sets of polynomials $\{P_{i}\}_{i\ge 0}$ and $\{Q_{j}\}_{j\ge 0}$,
 the so-called connection problem
between these polynomials is to determine the coefficients $A_{i,j}$ in the
expression
\[P_{i}(x)=\sum_{j=0}^{i}A_{i,j}\, Q_{j}(x).\]
The connection coefficients $A_{i,j}$ play an important role in many problems in pure and applied
mathematics and in mathematical physics. The problem of finding connection
coefficients between two sets of orthogonal polynomials has been
investigated by many authors, see for instance
\cite{area2001solving,doha2003connection,doha2004construction,doha2004recurrences,maroni2008connection,sanchez2001some,szwarc1992linearization}.
In fact,  most of the formulae for the connection coefficients
between orthogonal polynomials are given in terms of terminating
hypergeometric series of various types. For example, the connection
formula of Jacobi-Jacobi polynomials is given in terms
 of a terminating hypergeometric series of the type $_3F_{2}(1)$, see \cite{doha2004construction}.
 
 In this article, we solve the connection problem between Fibonacci polynomials and
 the orthogonal polynomials $T_n(x)$ and $U_n(x)$. We exhibit the Chebyshev expression form of Fibonacci polynomials $F_n(x)$. In fact, the connection coefficients turn out to be terminating hypergeometric series  of type ${_2F_1}(\lambda)$ where $\lambda$ is either $-4$ or $-1/4$. Furthermore, we tackle the inverse connection coefficients problem. In other words, we express $T_n(x)$ and $U_n(x)$ in terms of $F_n(x)$. Again, the latter coefficients involve the terminating hypergeometric series ${_2F_1}(\lambda)$.

Based on the new connection formulae that we derive, we display several identities satisfied by Fibonacci numbers, and we evaluate some finite sums. Moreover, we find relations between terminating hypergeometric series of type ${_2F_1}(\lambda)$ for certain values of $\lambda$ and specific parameters. Also, some
identities involving the derivatives sequences of
Fibonacci numbers are given.  
Finally, we express weighted definite integrals of products of Fibonacci and Chebyshev polynomials, and products of Fibonacci polynomials as sums involving hypergeometric series.  
 
\section{Some relevant properties of Fibonacci and Chebyshev polynomials}
 In this section, we recall the properties of Fibonacci and Chebyshev
 polynomials that we are going to use throughout the note.
  \subsection{Fibonacci polynomials}
 Fibonacci polynomials are generated via the following recurrence relation
 \[F_{n+2}(x)=x\, F_{n+1}(x)+F_{n}(x),\, n\ge 0,\textrm{ where }F_{0}(x)=0,\quad F_{1}(x)=1.\] 
   The $n$-th Fibonacci polynomial can be described explicitly as follows
    \[F_{n}(x)=\frac{\left(x+\sqrt{x^2+4}\right)^n-
    \left(x-\sqrt{x^2+4}\right)^n}{2^n\, \sqrt{x^2+4}}.
    \]
 The latter expression has the following explicit power form representation
 \begin{eqnarray*}
 F_{n}(x)&=&\sum_{j=0}^{\left\lfloor
 \frac{n-1}{2}\right\rfloor}\binom{n-j-1}{j}\, x^{n-2j-1}
 \end{eqnarray*}
 where $\left\lfloor
z\right\rfloor$ represents the largest integer less than or
equal to $z$.
  Now one can define the $n$-th Fibonacci number as follows
  \[F_{n}=F_{n}(1)=\frac{(1+\sqrt 5)^n-(1-\sqrt 5)^n}{2^n\, \sqrt 5}.\]
The corresponding derivatives sequences of the Fibonacci numbers are
denoted by $F^{(q)}_{n}$. They are defined via
   \[F^{(q)}_{n}=D^q\, F_{n}(x)\big|_{x=1}.\]
    Some of the identities relating Fibonacci and Lucas numbers to complex values of Chebyshev polynomials of the first and second
      kinds are as follows, see \cite{burcu2011complex,zhang2002chebyshev},
       \begin{align}
        \label{complex1}
        F_{n+1}&=\frac{1}{i^n}U_{n}\left(\frac{i}{2}\right),\\
      \label{complex2}
        U_{n}\left(-2\, i\right)&=\frac{(-i)^n}{2}\, F_{3(n+1)}.
         \end{align}
   There are several articles
   that study relations and
   identities satisfied by Fibonacci numbers and their derivatives sequences, see for example
   \cite{falcon2007fibonacci,dilcher2000hypergeometric,koshy2011fibonacci}.
    \subsection{Chebyshev polynomials of the first and second kinds}
    Chebyshev polynomials $T_{n}(x)$ and $U_{n}(x)$ of
 the first and second kinds, respectively, are polynomials in $x$, which can be
 defined by (see, Mason and Handscomb \cite{mason2010chebyshev}):
\[
 T_{n}(x)=\cos (n\ \theta),
\]
 and
\[
 U_{n}(x)=\frac{\sin (n+1)\theta}{\sin \theta},\]
 where $x=\cos\theta$.
 The polynomials\ $T_{n}(x)$ and $U_{n}(x)$ are orthogonal on $(-1,1)$
with respect to the weight functions $\frac{1}{\sqrt{1-x^2}}$ and
$\sqrt{1-x^2}$, that is
  \begin{equation} \label{ortho1}
 \int_{-1}^{1}
  \frac{1}{\sqrt{1-x^2}}\ T_{n}(x)\, T_{m}(x)\, dx
=\begin{cases}0, &m\not=n,\\
 \frac{\pi}{2}\, c_{n}, & m=n,\end{cases}
 \end{equation}
  and
  \begin{equation} \label{ortho2}
\int_{-1}^{1}\sqrt{1-x^2}\ U_{n}(x)\, U_{m}(x)\, dx
=\begin{cases}0, &m\not=n,\\
 \frac{\pi}{2}\, & m=n,\end{cases}
 \eq
  where
\[c_{n}=
 \begin{cases}
 2&n=0,\\
 1,&n>0.
 \end{cases}\]
  The polynomials $T_{n}(x)$ and $U_{n}(x)$ may be
generated, respectively, by means of the two recurrence relations:
\[
 T_{n}(x)=2x\, T_{n-1}(x)-T_{n-2}(x),\quad
n=2,3,\dots,
\]
   with
 \[T_{0}(x)=1,\quad T_{1}(x)=x,\]
 and
\[ U_{n}(x)=2x\, U_{n-1}(x)-U_{n-2}(x),\quad
n=2,3,\dots,\]
 with
  \[U_{0}(x)=1,\quad U_{1}(x)=2x.\]
   They also have the following explicit power forms:
\[T_{n}(x)=\frac{n}{2} \sum_{r=0}^{\left\lfloor
 \frac{n}{2}\right\rfloor}\frac{(-1)^r}{n-r}\binom{n-r}{r}\, (2\, x)^{n-2\, r},\]
  and
  \[U_{n}(x)=\sum_{r=0}^{\left\lfloor
 \frac{n}{2}\right\rfloor}(-1)^r\binom{n-r}{r}\, (2\, x)^{n-2\, r},\]
   The following special values are of important use later:
   \bq
   \label{inV}
  T_n(1)=1,\quad U_{n}(1)=n+1,
  \eq
  \bq
   \label{DerT}
  D^q
   T_{n}(1)=\prod_{i=0}^{q-1}\frac{(n-i)(n+i)}{2\, i+1},\
q\ge 1,
  \eq
  and
   \bq
   \label{DerU}
  D^q
   U_{n}(1)=(n+1)\prod_{i=0}^{q-1}\frac{(n-i)(n+i+2)}{2\, i+3},\
q\ge 1.
  \eq

\section{New connection formulae between
 Chebyshev polynomials of the first and second kinds and Fibonacci polynomials}
 The following two theorems establish two new connection
 formulae between Fibonacci polynomials and Chebyshev polynomials of
 the first and second kinds. The formulae are given in terms of values of
 hypergeometric functions.
\begin{Theorem}
\label{Thm1}
 For every $j\ge 1$, the following connection formula
holds:
 \begin{eqnarray*}
  \label{relation1}
  T_{j}(x)=j\, \sum_{m=0}^{\left\lfloor \frac{j}{2}\right\rfloor}
  (-1)^m\  \binom{j-m}{j-2m}\frac{2^{j-2 m-1}}{j-m}\
  \! _{2} F_{1}\left.\left(
\begin{array}{cc}
 -m, j-m\\ j-2m+2
\end{array}
\right|-4\right)\ F_{j-2m+1}(x).
 \end{eqnarray*}
 \end{Theorem}
    \begin{Theorem}
 \label{Thm2}
  For every $j\ge 1$, the following connection formula
holds:
 \begin{eqnarray*}
  \label{relation2}
  U_{j}(x)=2^j\, \sum_{m=0}^{\left\lfloor \frac{j}{2}\right\rfloor}
  (-1)^{m+1}\  \binom{j}{m}\ \frac{-j+2 m-1}{j-m+1}\
  \! _{2} F_{1}\left.\left(\begin{array}{cc}
 -m, -j+m-1\\ -j
\end{array}
\right|\frac{-1}{4}\right)\ F_{j-2m+1}(x).
 \end{eqnarray*}
\end{Theorem}

We will prove Theorem \ref{Thm1}. The proof of Theorem \ref{Thm2} is similar. 

\begin{Proof}
  The connection coefficients are written using the following terminating hypergeometric series  
\[_{2}F_{1}\left.\left(
\begin{array}{cc}
 -m, j-m\\ 2+j-2m
\end{array}
\right|-4\right)=\sum_{k=0}^m\frac{(-m)_k\,(j-m)_k\,
(-4)^k}{(2+j-2m)_k\,k!},\]
  Therefore it suffices to show that the following identity holds
  \begin{eqnarray*}
 \label{thm1eq1}
   T_j(x)=\phi_j(x):=\sum_{m=0}^{\left\lfloor
\frac{j}{2}\right\rfloor}
\sum_{k=0}^m\frac{(-1)^m\,2^{-1+j+2k-2m}\,j\,(1+j-2m)\,(j+k-m-1)!}
{k!\,(m-k)!\,(1+j+k-2m)!\,(m-k)!}F_{j-2m+1}(x).
 \end{eqnarray*}
 In other words, we want to prove that the function $\phi_j(x)$ satisfies the recurrence relation defining $T_j(x)$. 
It is easy to see that $\phi_{1}(x)=x$ and $\phi_{2}(x)=2x^2-1$.
We now make   use of the recurrence relation
  \[x\,F_j(x)=F_{j+1}(x)-F_{j-1}(x)\] satisfied by Fibonacci polynomials
  \begin{eqnarray*}
\label{thm1eq4}
 2x\,\phi_j(x)&=&\sum_{m=0}^{\left\lfloor
\frac{j}{2}\right\rfloor}
\sum_{k=0}^m\frac{(-1)^m\,2^{j+2k-2m}\,j\,(1+j-2m)\,(j+k-m-1)!}
{k!\,(m-k)!\,(1+j+k-2m)!\,(m-k)!}F_{j-2m+2}(x)\nonumber\\
 &+&\sum_{m=0}^{\left\lfloor
\frac{j}{2}\right\rfloor}
\sum_{k=0}^m\frac{(-1)^{m+1}\,2^{j+2k-2m}\,j\,(1+j-2m)\,(j+k-m-1)!}
{k!\,(m-k)!\,(1+j+k-2m)!\,(m-k)!}F_{j-2m}(x).
  \end{eqnarray*}
  Some algebraic manipulations will yield 
\begin{eqnarray*} 2x\,\phi_j(x)&=&\sum_{m=0}^{\left\lfloor
\frac{j-1}{2}\right\rfloor}
\sum_{k=0}^m\frac{(-1)^m\,2^{-2+j+2k-2m}\,(j-1)\,(j-2m)\,(j+k-m-2)!}
{k!\,(m-k)!\,(j+k-2m)!\,(m-k)!}F_{j-2m}(x)\\&+&\sum_{m=0}^{\left\lfloor
\frac{j+1}{2}\right\rfloor}
\sum_{k=0}^m\frac{(-1)^m\,2^{j+2k-2m}\,(j+1)\,(2+j-2m)\,(j+k-m)!}
{k!\,(m-k)!\,(2+j+k-2m)!\,(m-k)!}F_{j-2m+2}(x)\\&=&\phi_{j-1}(x)
+\phi_{j+1}(x),\end{eqnarray*}
 Since Chebyshev polynomials $T_{j}(x),\, j\ge 1$, are uniquely determined by the
 recurrence relation
  \[2\, x\, T_{j}(x)=T_{j-1}(x)+T_{j+1}(x),\, T_{1}=x,\, T_{2}(x)=2\, x^2-1,\] it follows that
  $\phi_{j}(x)$ is the $j$-th Chebyshev polynomial $ T_{j}(x)$ for $j\ge 1$. This completes the proof of Theorem \ref{Thm1}.
\end{Proof}
\section{Inversion formulae between Fibonacci and Chebyshev polynomials
   of the first and second kinds}
In this section we are concerned with deriving the inversion formulae to
those given in the previous section. We will express Fibonacci and Lucas polynomials in terms of Chebyshev polynomials. Again the connection coefficients turn out to be expressions involving values of hypergeometric functions of
the type $_2F_{1}$. We will need the following lemma in order to proceed.

 \begin{Lemma}
  \label{lemRec}
 We set $ d_{j,m}=\!{_2}F_{1}\left.\left(\begin{array}{cc}
 -m, -j+m-1\\ -j
\end{array}
\right|-4\right)$. The following recurrence relation
holds
 \begin{align*}
  \label{Rec}
4\ \binom{j-2}{m-1}\   (j-2 m+1)\  (j-m+1)\  d_{j-2,m-1}+
\binom{j-1}{m-1} \ (j-2 m+2)\  (j-m)\ d_{j-1,m-1}\\
+\binom{j-1}{m}\ (j-2 m)\  (j-m+1)
 d_{j-1,m}-\binom{j}{m} \ (j-m)\ (j-2 m+1)\ d_{j,m}=0.
  \end{align*}
  \end{Lemma}

\begin{Proof}
Let
 \[e_{j,m}=\binom{j}{m}\, (j-2m+1) \,d_{j,m}. \] 
 In order to prove that the recurrence relation above is satisfied, it suffices to show that
 \begin{eqnarray*}
  \label{eqRec}
4(j-m+1)\, e_{j-2,m-1}+(j-m)\, e_{j-1,m-1}+(j-m+1)\,
e_{j-1,m}=(j-m)\, e_{j,m}.
  \end{eqnarray*}
Now, $e_{j,m},$ can be written in the form
\[e_{j,m}=\binom{j}{m}\, (j-2m+1)\sum_{k=0}^m
\frac{(-m)_k\,(m-j-1)_k\, (-4)^k}{(-j)_k\ k!}.\]
  Using the identity
\[(-m)_k=\frac{(-1)^k\,m!}{(m-k)!},\]
$e_{j,m}$ can be written equivalently as
\[e_{j,m}=(j-m+1)(j-2m+1)\,\sum_{k=0}^m
\frac{(j-k)!\, 4^k}{k!\,(m-k)!\,(j-k-m+1)!}.\]
 It follows that 
  {\footnotesize\begin{eqnarray*}
   \label{steps1}
 4(j-m+1)\, e_{j-2,m-1}+(j-m)\, e_{j-1,m-1}+(j-m+1)\,
e_{j-1,m}&=&\\
(j-m)_2\bigg[4(1+j-2m) \sum_{k=0}^{m-1}\frac{(j-k)!\,
4^k}{k!\,(m-k)!\,(j-k-m)!}&+&(j-2m)\sum_{k=0}^{m}\frac{(j-k-1)!\,
4^k}{k!\,(m-k)!\, (j-k-m)!}
\\ \qquad \qquad +(2+j-2m)\sum_{k=0}^{m-1}
\frac{(j-k-1)!\, 4^k}{k!\,(m-k-1)!\,(j-k-m+1)!} \bigg].
  \end{eqnarray*}}
 Taking a common denominator yields the following simplification
 \begin{eqnarray*}
 4(j-m+1)\, e_{j-2,m-1}+(j-m)\, e_{j-1,m-1}+(j-m+1)\,
e_{j-1,m}&=&\\
(j-m)_2\,
(1+j-2m)\sum_{k=0}^{m}\frac{(j-k)!}{k!\,(m-k-1)!\,(j-k-m+1)!}
&=&(j-m)\, e_{j,m}.
 \end{eqnarray*}
 Lemma \ref{lemRec} is now proved.
 \end{Proof}
 
 In the following two theorems, we exhibit the
 connection relation between Fibonacci polynomials and Chebyshev
 polynomials of the first and second kinds.

  \begin{Theorem}
   \label{thm3}
    For every nonnegative integer $j$, the following connection
    formula holds
 \begin{eqnarray*}
  F_{j+1}(x)=\sum_{m=0}^{\left\lfloor
  \frac{j}{2}\right\rfloor}\frac{1}{c_{j-2m}}\ \binom{j-m}{j-2m}\ 2^{-j+2 m+1}\
  \! _{2}F_{1}\left.\left(\begin{array}{cc}
 -m, j-m+1\\ j-2m+1
\end{array}
\right|\frac{-1}{4}\right)\ T_{j-2m}(x),
 \end{eqnarray*}
  where
   \[ c_{j}=
 \begin{cases}
 2&j=0,\\
 1,&j>0.
 \end{cases}\]
\end{Theorem}
  \begin{Theorem}
   \label{thm4}
       For every nonnegative integer $j$, the following connection
    formula holds
 \begin{eqnarray*}
  F_{j+1}(x)=\frac{1}{2^j}\, \sum_{m=0}^{\left\lfloor \frac{j}{2}\right\rfloor}
  \binom{j}{m}\ \frac{(j-2 m+1)}{j-m+1}\
  \! _{2}F_{1}\left.\left(\begin{array}{cc}
 -m, -j+m-1\\ -j
\end{array}
\right|-4\right)\ U_{j-2m}(x).
 \end{eqnarray*}
\end{Theorem}

The proofs of Theorems \ref{thm3} and \ref{thm4} are similar, so it suffices to
 prove Theorem \ref{thm4}.

\begin{Proof}
  We will proceed by induction. Assume that the above identity is
  valid for any $k\le j$. 
  We know that 
 \[ F_{j+1}(x)=x\, F_{j}(x)+F_{j-1}(x),\]
 therefore using the induction hypothesis to write $F_{j-1}(x)$ and $F_j(x)$ in terms of Chebyshev polynomials, together with making use of
the recurrence relation
\[x\, U_{j}(x)=\frac12\left[U_{j-1}(x)+U_{j+1}(x)\right],\]
yield

{\footnotesize \begin{eqnarray*}
 \label{Indstep}
F_{j+1}(x)&=&\frac{1}{2^j}\sum_{m=0}^{\left\lfloor
\frac{j-1}{2}\right\rfloor}\binom{j-1}{m}\frac{j-2 m}{j-m} \
_{2}F_{1}\left.\left(\begin{array}{cc}
 -m, -j+m\\ 1-j
\end{array}
\right|-4\right) \big(U_{j-2 m-2}(x)+U_{j-2 m}(x)\big)\\&+&
 \frac{1}{2^{j-2}}\sum_{m=0}^{\left\lfloor
\frac{j}{2}\right\rfloor-1}\frac{j-2 m-1}{j-m-1}\
 \binom{j-2}{m} \, _{2}F_{1}\left.\left(\begin{array}{cc}
 -m, -j+m+1\\ 2-j
\end{array}
\right|-4\right) U_{j-2 m-2}(x).
 \end{eqnarray*}}
  In fact the latter identity can be simplified and rewritten as follows
\begin{eqnarray*}
 \label{Simplify1}
F_{j+1}(x)=\sum_{m=0}^{\left\lfloor
\frac{j-1}{2}\right\rfloor}g_{j,m}U_{j-2m}(x)+\frac12\,
g_{j-1,\mu_{j}}\, U_{j-2\, \mu_{j}-2}(x)+g_{j-2,\nu_{j}-1}\
U_{j-2\nu_{j}}(x)\ \theta_{j},
 \end{eqnarray*}
where
 {\footnotesize\begin{eqnarray*}
 g_{j,m}=\frac{1}{2^j}\left\{\binom{j-1}{m-1}\ \frac{(j-2 m+2)}
 {j-m+1}\ _{2}F_{1}\left.\left(\begin{array}{cc}
 1-m,-j+m-1\\ 1-j
\end{array}
\right|-4\right)\right.+\\
  4\, \binom{j-2}{m-1}\ \frac{(j-2 m+1)}
 {j-m}\ _{2}F_{1}\left.\left(\begin{array}{cc}
 1-m,m-j\\ 2-j
\end{array}
\right|-4\right)+\\
  \left.\binom{j-1}{m}\ \frac{(j-2 m)}
 {j-m}\ _{2}F_{1}\left.\left(\begin{array}{cc}
 -m,m-j\\ 1-j
\end{array}
\right|-4\right)\right\},
 \end{eqnarray*}}
  and
 \[\mu_{j}=\left\lfloor
\frac{j-1}{2}\right\rfloor,\ \nu_{j}=\left\lfloor
\frac{j}{2}\right\rfloor,\ \theta_{j}=\left\{\begin{array}{ll}1,
\textrm{ if j\  even},\\0,   \textrm{ if j\  odd}.\end{array}\right.\]

Making use of Lemma \ref{lemRec} together with some manipulations imply that $g_{j,m}$ can be reduced to take the
form
\[g_{j,m}=\frac{\binom{j}{m}\, (j-2m+1)}{2^j (j-m+1)}\ _{2}F_{1}
\left.\left(\begin{array}{cc}
 -m,-1+m-j\\ -j
\end{array}
\right|-4\right).\]
 This yields that 
\[
  \begin{split}
  F_{j+1}(x)=\frac{1}{2^j}\, \sum_{m=0}^{\left\lfloor \frac{j}{2}\right\rfloor}
  \frac{\binom{j}{m}\ (j-2 m+1)}{j-m+1}\
  \! _{2}F_{1}\left.\left(\begin{array}{cc}
 -m, -j+m-1\\ -j
\end{array}
\right|-4\right)\ U_{j-2m}(x),
\end{split}\]
 which completes the proof of Theorem \ref{thm4}.
 \end{Proof}

  \section {Some applications}
    In this section, we introduce three applications to the new
    derived connection formulae and their inversion ones: (i) We display some new expressions involving Fibonacci and Lucas
    numbers. (ii) Some new expressions for the derivatives sequences of
    Fibonacci numbers are given. (iii) We
    evaluate some definite integrals involving certain
   products of Fibonacci and Chebyshev polynomials.
   
\subsection{New expressions involving Fibonacci and Lucas numbers}
 In this section we use the results we developed 
 in \S 3 and \S 4 to evaluate finite sums involving certain values of hypergeometric functions and Fibonacci numbers.
  \begin{Corollary}
   \label{cor1}
    For every nonnegative integer $j$, the following two identities
    hold:
\begin{eqnarray*}
  \label{num1}
  \begin{split}
  \sum_{m=0}^{\left\lfloor \frac{j}{2}\right\rfloor}
  (-1)^m\  \binom{j-m}{j-2m}\frac{2^{j-2 m-1}}{j-m}\
  \! _{2} F_{1}\left.\left(
\begin{array}{cc}
 -m, j-m\\ j-2m+2
\end{array}
\right|-4\right)\ F_{j-2m+1}=1,
\end{split}
 \end{eqnarray*}
  and
\begin{eqnarray*}
  \label{num2}
  \begin{split}
  2^j\, \sum_{m=0}^{\left\lfloor \frac{j}{2}\right\rfloor}
 (-1)^{m+1}\  \binom{j}{m}\  \frac{(-j+2 m-1)}{j-m+1}\
  \! _{2} F_{1}\left.\left(\begin{array}{cc}
 -m, -j+m-1\\ -j
\end{array}
\right|\frac{-1}{4}\right)\ F_{j-2m+1}=j+1.
\end{split}
 \end{eqnarray*}
  \end{Corollary}
   \begin{Corollary}
    \label{cor2}
     For every nonnegative integer $j$, the following two expressions for Fibonacci
     numbers are valid
 \begin{eqnarray*}
  \label{num3}
  \begin{split}
  F_{j+1}=\sum_{m=0}^{\left\lfloor
  \frac{j}{2}\right\rfloor}\frac{1}{c_{j-2m}}\ \binom{j-m}{j-2m}\ 2^{-j+2 m+1}\
  \! _{2}F_{1}\left.\left(\begin{array}{cc}
 -m, j-m+1\\ j-2m+1
\end{array}
\right|\frac{-1}{4}\right),
\end{split}
 \end{eqnarray*}

 \begin{eqnarray*}
  \label{num4}
  \begin{split}
  F_{j+1}=\frac{1}{2^j}\, \sum_{m=0}^{\left\lfloor
  \frac{j}{2}\right\rfloor} \binom{j}{m}\
  \frac{(j-2m+1)^2}{j-m+1}\
  \! _{2}F_{1}\left.\left(\begin{array}{cc}
 -m, -j+m-1\\ -j
\end{array}
\right|-4\right).
\end{split}
 \end{eqnarray*}
  \end{Corollary}
   \begin{Proof}
    The proof of Corollaries \ref{cor1} and \ref{cor2} are immediately obtained
    from Theorems \ref{Thm1}, \ref{Thm2}, \ref{thm3} and
    \ref{thm4}, respectively, by setting $x=1$.
   \end{Proof}
    The fact that Fibonacci numbers themselves can be expressed as values of terminating hypergeometric series can be exploited in order to find a linear relation between hypergeometric series of type ${_2 F_1}$ with different parameters at a specific value.
    
    The $n$-th Fibonacci number can be written as a hypergeometric series itself. In fact one knows that \begin{eqnarray*} F_n&=&\frac{n}{2^{n-1}}\;{_2 F_1}\left.\left(\begin{array}{cc}
 \frac{1-n}{2}, \frac{2-n}{2}\\ \frac{3}{2}
\end{array}
\right|5\right)
={_2 F_1}\left.\left(\begin{array}{cc}
 \frac{1-n}{2}, \frac{2-n}{2}\\ 1-n
\end{array}
\right|-4\right),
\end{eqnarray*} see \cite{dilcher2000hypergeometric} for example. One of the linear transformations of hypergeometric series is given by the following identity
\[{_2 F_1}\left.\left(\begin{array}{cc}
 a, b\\ c
\end{array}
\right|z\right)=(1-z)^{-a}\;{_2 F_1}\left.\left(\begin{array}{cc}
 a, c-b\\ c
\end{array}
\right|\frac{z}{z-1}\right).\] Putting these together, one can rewrite Corollary \ref{cor2} as follows. 
    
 \begin{Corollary}
    \label{cor22}
     For every nonnegative integer $j$, the following two expressions identities hold true
     \begin{eqnarray*}
  {_2 F_1}\left.\left(\begin{array}{cc}
 \frac{-j}{2}, \frac{1-j}{2}\\ -j
\end{array}
\right|-4\right)&=&\sum_{m=0}^{\left\lfloor
  \frac{j}{2}\right\rfloor}\frac{1}{c_{j-2m}}\ \binom{j-m}{j-2m}\ 2^{-j+2 m+1}\
  \! _{2}F_{1}\left.\left(\begin{array}{cc}
 -m, j-m+1\\ j-2m+1
\end{array}
\right|\frac{-1}{4}\right)\\
&=&
\frac{1}{2^j}\, \sum_{m=0}^{\left\lfloor
  \frac{j}{2}\right\rfloor} \binom{j}{m}\
  \frac{(j-2m+1)^2}{j-m+1}\
  \! _{2}F_{1}\left.\left(\begin{array}{cc}
 -m, -j+m-1\\ -j
\end{array}
\right|-4\right),\\
  {_2 F_1}\left.\left(\begin{array}{cc}
 \frac{-j}{2}, \frac{1-j}{2}\\ \frac{3}{2}
\end{array}
\right|5\right)&=&\frac{2}{j+1}\sum_{m=0}^{\left\lfloor
  \frac{j}{2}\right\rfloor}\frac{5^m}{c_{j-2m}}\ \binom{j-m}{j-2m}\
  \! _{2}F_{1}\left.\left(\begin{array}{cc}
 -m, -m\\ j-2m+1
\end{array}
\right|\frac{1}{5}\right),\\
&=&\frac{1}{2^j}\, \sum_{m=0}^{\left\lfloor
  \frac{j}{2}\right\rfloor} 5^m\binom{j}{m}\
  \frac{(j-2m+1)^2}{j-m+1}\
  \! _{2}F_{1}\left.\left(\begin{array}{cc}
 -m, 1-m\\ -j
\end{array}
\right|\frac{4}{5}\right).
 \end{eqnarray*}
 
  \end{Corollary}
Furthermore, one may obtain similar identities to the ones above relating values of the hypergeometric series ${_2 F_1}$ evaluated at $1/5$ or $4/5$ and several other values, see \S 4 in \cite{dilcher2000hypergeometric}, using different linear transformations.
    
     Now, and based on the identities \eqref{complex1} and \eqref{complex2},
     along with the connection formulae obtained in \S 4, the following identities follow.
     \begin{Corollary}
      For every nonnegative integer $j$, and $i^2=-1$, the following
      identities hold:
      {\footnotesize \begin{eqnarray*}
     i^j\, F_{j+1}&=&2^{j}\, \sum_{m=0}^{\frac{j}{2}}
     \frac{(-1)^{m+1}\ \binom{j}{m}  (-j+2 m-1) F_{j-2 m+1}
     \left(\frac{i}{2}\right)}{j-m+1}\
      _{2} F_{1}\left.\left(\begin{array}{cc}
 -m, -j+m-1\\ -j
\end{array}
\right|\frac{-1}{4}\right),\\
     (-i)^j\, F_{3\, j+3}&=&2^{j+1}\, \sum_{m=0}^{\frac{j}{2}}
     \frac{(-1)^{m+1}\binom{j}{m}\  (-j+2 m-1) \ F_{j-2 m+1}(-2
     i)}{j-m+1}\
      _{2} F_{1}\left.\left(\begin{array}{cc}
 -m, -j+m-1\\ -j
\end{array}
\right|\frac{-1}{4}\right).
\end{eqnarray*}}
\end{Corollary}
One may also obtain some new trigonometric identities making use of the fact that $T_n(\cos \theta)=\cos(n\theta)$. In fact one has 
\begin{eqnarray*}
  F_{j+1}(\cos\theta)=\sum_{m=0}^{\left\lfloor
  \frac{j}{2}\right\rfloor}\frac{1}{c_{j-2m}}\ \binom{j-m}{j-2m}\ 2^{-j+2 m+1}\
  \! _{2}F_{1}\left.\left(\begin{array}{cc}
 -m, j-m+1\\ j-2m+1
\end{array}
\right|\frac{-1}{4}\right)\ \cos\left((j-2m)\theta\right),
 \end{eqnarray*}
  where
   \[ c_{j}=
 \begin{cases}
 2&j=0,\\
 1,&j>0.
 \end{cases}.\]
 Another interesting identity is obtained by observing that $\displaystyle T_n\left(\frac{x+x^{-1}}{2}\right)=\frac{x^n+x^{-n}}{2}$, whence 
 {\footnotesize\begin{eqnarray*}
  F_{j+1}\left(\frac{x+x^{-1}}{2}\right)=\sum_{m=0}^{\left\lfloor
  \frac{j}{2}\right\rfloor}\frac{1}{c_{j-2m}}\ \binom{j-m}{j-2m}\ 2^{-j+2 m}\
  \! _{2}F_{1}\left.\left(\begin{array}{cc}
 -m, j-m+1\\ j-2m+1
\end{array}
\right|\frac{-1}{4}\right)\ \left( x^{j-2m}+x^{2m-j}\right).
 \end{eqnarray*}}
 
   \subsection{New derivatives sequences identities}
    Based on the connection formulae introduced in Theorems \ref{Thm1}, \ref{Thm2}, \ref{thm3}, and \ref{thm4}, we may obtain new formulae for the derivatives sequences of
    Fibonacci numbers by using the identities \eqref{DerT} and \eqref{DerU}.
     \begin{Corollary}
     \label{corder1}
      For all $q\ge 1$, the following two formulae are valid
    \begin{eqnarray*}
    \label{Derq1}
    \begin{split}
&\sum_{m=0}^{\left\lfloor \frac{j}{2}\right\rfloor}(-1)^m\,
\binom{j-m}{j-2 m}\  \frac{ 2^{j-2 m-1} }{j-m}\ _2 F_{1}\left.\left(
\begin{array}{cc}
 -m, j-m\\ j-2m+2
\end{array}
\right|-4\right)\, F^{(q)}_{j-2m+1}=\\
&\hspace{40pt}\frac{\ (-1)^{q+1} \ \sqrt{\pi }\  j\  (1-j)_{q-1}
(j+1)_{q-1}}{2^q\ \Gamma \left(q+\frac{1}{2}\right)},
\end{split}
 \end{eqnarray*}
  and
 \begin{eqnarray*}
  \label{Derq2}
  \begin{split}
  &\sum_{m=0}^{\left\lfloor \frac{j}{2}\right\rfloor}
 (-1)^{m+1}\  \binom{j}{m}\  \frac{(-j+2 m-1)}{j-m+1}\
  \! _{2}F_{1}\left.\left(\begin{array}{cc}
 -m, -j+m-1\\ -j
\end{array}
\right|\frac{-1}{4}\right)\, F^{(q)}_{j-2m+1}=\\
&\hspace{40pt}\frac{\  (-1)^{q+1}\sqrt{\pi }\  (j)_{3}\, (1-j)_{q-1}
(j+3)_{q-1}}{2^{j+q+1} \ \Gamma \left(q+\frac{3}{2}\right)}.
\end{split}
 \end{eqnarray*}
  \end{Corollary}
  \begin{Corollary}
  \label{corder2}
      For all $q\ge 1$, the following two formulae are valid
 \begin{eqnarray*}
  \label{Derq3}
  \begin{split}
  F^{(q)}_{j+1}=\frac{(-1)^{q+1}\ \sqrt{\pi}}{\Gamma \left(q+\frac{1}{2}\right)}\
  &\sum_{m=0}^{\left\lfloor \frac{j}{2}\right\rfloor}\frac{1}{c_{j-2m}}\
  \binom{j-m}{j-2m}\ 2^{-j+2 m-q+1}\ (j-2 m)^2 \ (j-2 m+1)_{q-1}\
    \times\\
  &(-j+2 m+1)_{q-1}\  \ _{2}F_{1}\left.\left(\begin{array}{cc}
 -m, j-m+1\\ j-2m+1
\end{array}
\right|\frac{-1}{4}\right),
\end{split}
 \end{eqnarray*}
  and
 \begin{eqnarray*}
  \label{Derq4}
  \begin{split}
  F^{(q)}_{j+1}=\frac{(-1)^{q+1}\ \sqrt{\pi }}{2^{j+q+1}\ \Gamma \left(q+\frac{3}{2}\right)}
  &\, \sum_{m=0}^{\left\lfloor \frac{j}{2}\right\rfloor}
 \binom{j}{m} \  \frac{(j-2 m) (j-2 m+1)^2\  (j-2 m+2)}{j-m+1}\times\\
  &(j-2
 m+3)_{q-1}\ \  _{2}F_{1}\left.\left(\begin{array}{cc}
 -m, -j+m-1\\ -j
\end{array}
\right|-4\right).
\end{split}
 \end{eqnarray*}
  \end{Corollary}
   \begin{Proof}
    In order to prove Corollaries \ref{corder1} and \ref{corder2},
  one needs to differentiate the
connection formulae in Theorems \ref{Thm1}, \ref{Thm2}, \ref{thm3}, and \ref{thm4}, then one sets $x=1$. Now the identities follow from \eqref{DerT} and \eqref{DerU}.
   \end{Proof}
 \subsection{Some integrals formulae involving Chebyshev and Fibonacci polynomials}
  The following two integrals formulae are direct consequences of
    Theorems \ref{thm3} and \ref{thm4}.
  \begin{Corollary}
     For all $j\ge k$, the following two integrals
   formulae hold:
  {\footnotesize\begin{eqnarray*}
  \label{int1}
  \int\limits_{-1}^{1}\frac{F_{j+1}(x) T_k(x)}{\sqrt{1-x^2}}\, dx=
      \begin{cases}
  \frac{\pi\ \binom{\frac{j+k}{2}}{k}}
  {2^k\, c_{k}}\ _{2}F_{1}\left.\left(\begin{array}{cc}
 \frac{k-j}{2},\frac{1}{2} (j+k+2)\\ k+1
\end{array}
\right|\frac{-1}{4}\right),& \textrm{if }(k+j)\ \text{even},\\
 0,& \text{otherwise},
 \end{cases}
  \end{eqnarray*}}
   and
     {\footnotesize\begin{eqnarray*}
  \label{int2}
  \int\limits_{-1}^{1}\sqrt{1-x^2}\ F_{j+1}(x)\  U_k(x)\, dx=
      \begin{cases}
  \frac{\pi\
  \binom{j}{\frac{j-k}{2}}\ (k+1)}{2^j(j+k+2)}\
  \ \! _{2}F_{1}\left.\left(\begin{array}{cc}
 \frac{k-j}{2},\frac{-1}{2} (j+k+2)\\ -j
\end{array}
\right|-4\right),& \textrm{if }(k+j)\ \text{even},\\
 0,& \text{otherwise}.
 \end{cases}
  \end{eqnarray*}}
\end{Corollary}
 Now one can find an explicit description for weighted definite integrals of products of Fibonacci polynomials in terms of hypergeometric series. 
 \begin{Corollary}
 For all $j\ge k$, the following two integrals
   formulae hold:
  {\footnotesize\begin{align*}
  \int\limits_{-1}^{1}\frac{F_{j+1}(x) F_{k+1}(x)}{\sqrt{1-x^2}}\, dx=&\\ \frac{\pi}{2^{k+j-1}}&\sum_{m=0}^{\lfloor\frac{k}{2}\rfloor}
  2^{4m}d_m\frac{{j-m\choose j-2m}{k-m\choose k-2m}}{c_{k-2m}c_{j-2m}} \;{_{2}F_{1}}\left.\left(\begin{array}{cc}
 -m,k-m+1\\ k-2m+1
\end{array}
\right|\frac{-1}{4}\right) \times{_{2}F_{1}}\left.\left(\begin{array}{cc}
 -m,j-m+1\\ j-2m+1
\end{array}
\right|\frac{-1}{4}\right)
  \end{align*}}
   and
  {\footnotesize\begin{align*}
  \int\limits_{-1}^{1}\sqrt{1-x^2}F_{j+1}(x) F_{k+1}(x)\, dx=&\\ \frac{\pi}{2^{k+j+1}}\sum_{m=0}^{\lfloor\frac{k}{2}\rfloor}
  {j\choose m}{k\choose m}&\frac{(k-2m+1)(j-2m+1)}{(k-m+1)(j-m+1)} \;{_{2}F_{1}}\left.\left(\begin{array}{cc}
 -m,-k+m-1\\ -k
\end{array}
\right|-4\right)\\ &\times{_{2}F_{1}}\left.\left(\begin{array}{cc}
 -m,-j+m-1\\ -j
\end{array}
\right|-4\right).
  \end{align*}}
 \end{Corollary}

\bibliographystyle{plain}
\bibliography{FibCheb}

\begin{thebibliography}{10}

\bibitem{WITSF}
W.M. Abd-Elhameed.
\newblock New product and linearization formulae of \textsc{J}acobi polynomials
  of certain parameters.
\newblock {\em Integr. Transf. Spec. F.}, DOI:10.1080/10652469.2015.1029924.

\bibitem{Ramanujan}
W.M. Abd-Elhameed, E.H. Doha, and H.M. Ahmed.
\newblock Linearization formulae for certain \textsc{J}acobi polynomials.
\newblock {\em Ramanujan J., DOI: 10.1007/s11139-014-9668-2}.

\bibitem{area2001solving}
I.~Area, E.~Godoy, A.~Ronveaux, and A.~Zarzo.
\newblock Solving connection and linearization problems within the askey scheme
  and its $q-$analogue via inversion formulas.
\newblock {\em J. Comput. Appl. Math.}, 133(1):151--162, 2001.

\bibitem{burcu2011complex}
S.B. Bozkurt, F.~Y{\i}lmaz, and D.~Bozkurt.
\newblock On the complex factorization of the \textsc{L}ucas sequence.
\newblock {\em Appl. Math. Lett.}, 24(8):1317--1321, 2011.

\bibitem{dilcher2000hypergeometric}
K.~Dilcher and L.~Pisano.
\newblock Hypergeometric functions and \textsc{F}ibonacci numbers.
\newblock {\em Fibonacci Quarterly}, 38(4):342--362, 2000.

\bibitem{doha2003connection}
E.H. Doha.
\newblock On the connection coefficients and recurrence relations arising from
  expansions in series of \textsc{L}aguerre polynomials.
\newblock {\em J. Phys. A: Math. Gen.}, 36(20):5449--5462, 2003.

\bibitem{doha2004construction}
E.H. Doha.
\newblock On the construction of recurrence relations for the expansion and
  connection coefficients in series of \textsc{J}acobi polynomials.
\newblock {\em J. Phys. A: Math. Gen.}, 37(3):657, 2004.

\bibitem{dohanew}
E.H. Doha and W.M. Abd-Elhameed.
\newblock New linearization formulae for the products of \textsc{C}hebyshev
  polynomials of third and fourth kind.
\newblock {\em Rocky Mt. J. Math., To appear.}

\bibitem{doha2014integrals}
EH~Doha and WM~Abd-Elhameed.
\newblock Integrals of \textsc{C}hebyshev polynomials of third and fourth
  kinds: An application to solution of boundary value problems with polynomial
  coefficients.
\newblock {\em J. Contemp. Math. Anal.}, 49(6):296--308, 2014.

\bibitem{doha2004recurrences}
E.H. Doha and H.M. Ahmed.
\newblock Recurrences and explicit formulae for the expansion and connection
  coefficients in series of \textsc{B}essel polynomials.
\newblock {\em J. Phys. A: Math. Gen.}, 37(33):8045, 2004.

\bibitem{el2010new}
M.~El-Mikkawy and T.~Sogabe.
\newblock A new family of $k-$\textsc{F}ibonacci numbers.
\newblock {\em Appl. Math. Comput.}, 215(12):4456--4461, 2010.

\bibitem{falcon2007fibonacci}
S.~Falcon and A.~Plaza.
\newblock On the \textsc{F}ibonacci $k-$numbers.
\newblock {\em Chaos Soliton Fract.}, 32(5):1615--1624, 2007.

\bibitem{gulec2013new}
H.H. Gulec, N.~Taskara, and K.~Uslu.
\newblock A new approach to generalized \textsc{F}ibonacci and \textsc{L}ucas
  numbers with binomial coefficients.
\newblock {\em Appl. Math. Comput.}, 220:482--486, 2013.

\bibitem{koshy2011fibonacci}
T.~Koshy.
\newblock {\em \textsc{F}ibonacci and \textsc{L}ucas numbers with
  applications}, volume~51.
\newblock John Wiley \& Sons, 2011.

\bibitem{maroni2008connection}
P.~Maroni and Z.~da~Rocha.
\newblock Connection coefficients between orthogonal polynomials and the
  canonical sequence: an approach based on symbolic computation.
\newblock {\em Numer. Algor.}, 47(3):291--314, 2008.

\bibitem{mason2010chebyshev}
J.C. Mason and D.C. Handscomb.
\newblock {\em \textsc{C}hebyshev \textsc{p}olynomials}.
\newblock Chapman and Hall, New York, NY, CRC, Boca Raton, 2010.

\bibitem{sanchez2001some}
J.~S{\'a}nchez-Ruiz and J.S. Dehesa.
\newblock Some connection and linearization problems for polynomials in and
  beyond the askey scheme.
\newblock {\em J. Comput. Appl. Math.}, 133(1):579--591, 2001.

\bibitem{szwarc1992linearization}
R.~Szwarc.
\newblock Linearization and connection coefficients of orthogonal polynomials.
\newblock {\em Mh. Math.}, 113(4):319--329, 1992.

\bibitem{taskara2010properties}
N.~Taskara, K.~Uslu, and H.H. Gulec.
\newblock On the properties of \textsc{L}ucas numbers with binomial
  coefficients.
\newblock {\em Appl. Math. Lett.}, 23(1):68--72, 2010.

\bibitem{wloch2013some}
A.~W{\l}och.
\newblock Some identities for the generalized \textsc{F}ibonacci numbers and
  the generalized \textsc{L}ucas numbers.
\newblock {\em Appl. Math. Comput.}, 219(10):5564--5568, 2013.

\bibitem{yazlik2012note}
Yasin Y.~Yazlik and N.~Taskara.
\newblock A note on generalized-\textsc{H}oradam sequence.
\newblock {\em Comput. Math. Appl.}, 63(1):36--41, 2012.

\bibitem{zhang2002chebyshev}
W.~Zhang.
\newblock On \textsc{C}hebyshev polynomials and \textsc{F}ibonacci numbers.
\newblock {\em Fibonacci Quart.}, 40(5):424--428, 2002.

\end{thebibliography}
\end{document}